# ASYMPTOTIC EQUIVALENCE AND ADAPTIVE ESTIMATION FOR ROBUST NONPARAMETRIC REGRESSION

BY T. TONY CAI[1] AND HARRISON H. ZHOU[2]

*University of Pennsylvania and Yale University*

Asymptotic equivalence theory developed in the literature so far are only for bounded loss functions. This limits the potential applications of the theory because many commonly used loss functions in statistical inference are unbounded. In this paper we develop asymptotic equivalence results for robust nonparametric regression with unbounded loss functions. The results imply that all the Gaussian nonparametric regression procedures can be robustified in a unified way. A key step in our equivalence argument is to bin the data and then take the median of each bin.

The asymptotic equivalence results have significant practical implications. To illustrate the general principles of the equivalence argument we consider two important nonparametric inference problems: robust estimation of the regression function and the estimation of a quadratic functional. In both cases easily implementable procedures are constructed and are shown to enjoy simultaneously a high degree of robustness and adaptivity. Other problems such as construction of confidence sets and nonparametric hypothesis testing can be handled in a similar fashion.

**1. Introduction.** The main goal of the asymptotic equivalence theory is to approximate general statistical models by simple ones. If a complex model is asymptotically equivalent to a simple model, then all asymptotically optimal procedures can be carried over from the simple model to the complex one for *bounded* loss functions and the study of the complex model is then essentially simplified. Early work on asymptotic equivalence theory was focused on the parametric models and the equivalence is local. See Le Cam (1986).

Received August 2008; revised December 2008.
[1]Supported in part by NSF Grant DMS-06-04954.
[2]Supported in part by NSF Grant DMS-06-45676.

There have been important developments in the asymptotic equivalence theory for nonparametric models in the last decade or so. In particular, global asymptotic equivalence theory has been developed for nonparametric regression in Brown and Low (1996b) and Brown et al. (2002), nonparametric density estimation models in Nussbaum (1996) and Brown et al. (2004), generalized linear models in Grama and Nussbaum (1998), nonparametric autoregression in Milstein and Nussbaum (1989), diffusion models in Delattre and Hoffmann (2002) and Genon-Catalot, Laredo and Nussbaum (2002), GARCH model in Wang (2002) and Brown, Wang and Zhao (2003), and spectral density estimation in Golubev, Nussbaum and Zhou (2009).

So far all the asymptotic equivalence results developed in the literature are only for bounded loss functions. However, for many statistical applications, asymptotic equivalence under bounded losses is not sufficient because many commonly used loss functions in statistical inference such as squared error loss are unbounded. As commented by Johnstone (2002) on the asymptotic equivalence results: "Some cautions are in order when interpreting these results.... Meaningful error measures... may not translate into, say, squared error loss in the Gaussian sequence model."

In this paper we develop asymptotic equivalence results for robust nonparametric regression with an unknown symmetric error distribution for *unbounded* loss functions which include, for example, the commonly used squared error and integrated squared error losses. Consider the nonparametric regression model

$$(1) \qquad Y_i = f\left(\frac{i}{n}\right) + \xi_i, \qquad i = 1, \ldots, n,$$

where the errors $\xi_i$ are independent and identically distributed with some density $h$. The error density $h$ is assumed to be symmetric with median 0, but otherwise unknown. Note that for some heavy-tailed distributions such as Cauchy distribution the mean does not even exist. We thus do not assume the existence of the mean here. One is often interested in robustly estimating the regression function $f$ or some functionals of $f$. These problems have been well studied in the case of Gaussian errors. In the present paper we introduce a unified approach to turn the general nonparametric regression model (1) into a standard Gaussian regression model and then in principle any procedure for Gaussian nonparametric regression can be applied. More specifically, with properly chosen $T$ and $m$, we propose to divide the observations $Y_i$ into $T$ bins of size $m$ and then take the median $X_j$ of the observations in the $j$th bin for $j = 1, \ldots, T$. The asymptotic equivalence results developed in Section 2 show that under mild regularity conditions, for a wide collection of error distributions the experiment of observing the medians $\{X_j : j = 1, \ldots, T\}$ is in fact asymptotically equivalent to the standard

Gaussian nonparametric regression model

$$(2) \qquad Y_i = f\left(\frac{i}{T}\right) + \frac{1}{2h(0)\sqrt{m}} z_i, z_i \overset{\text{i.i.d.}}{\sim} N(0,1), \qquad i = 1, \ldots, T$$

for a large class of unbounded losses. Detailed arguments are given in Section 2.

We develop the asymptotic equivalence results for the general regression model (1) by first extending the classical formulation of asymptotic equivalence in Le Cam (1964) to accommodate unbounded losses. The asymptotic equivalence result has significant practical implications. It implies that all statistical procedures for any asymptotic decision problem in the setting of the Gaussian nonparametric regression can be carried over to solve problems in the general nonparametric regression model (1) for a class of unbounded loss functions. In other words, all the Gaussian nonparametric regression procedures can be robustified in a unified way. We illustrate the applications of the general principles in two important nonparametric inference problems under the model (1): robust estimation of the regression function $f$ under integrated squared error loss and the estimation of the quadratic functional $Q(f) = \int f^2$ under squared error.

As we demonstrate in Sections 3 and 4 the key step in the asymptotic equivalence theory, binning and taking the medians, can be used to construct simple and easily implementable procedures for estimating the regression function $f$ and the quadratic functional $\int f^2$. After obtaining the medians of the binned data, the general model (1) with an unknown symmetric error distribution is turned into a familiar Gaussian regression model, and then a Gaussian nonparametric regression procedure can be applied. In Section 3 we choose to employ a blockwise James–Stein wavelet estimator, BlockJS, for the Gaussian regression problem because of its desirable theoretical and numerical properties. See Cai (1999). The robust wavelet regression procedure has two main steps:

1. Binning and taking median of the bins.
2. Applying the BlockJS procedure to the medians.

The procedure is shown to achieve four objectives simultaneously: robustness, global adaptivity, spatial adaptivity, and computational efficiency. Theoretical results in Section 3.2 show that the estimator achieves optimal global adaptation for a wide range of Besov balls as well as a large collection of error distributions. In addition, it attains the local adaptive minimax rate for estimating functions at a point. Figure 1 compares a direct wavelet estimate with our robust estimate in the case of Cauchy noise. The example illustrates the fact that direct application of a wavelet regression procedure designed for Gaussian noise may not work at all when the noise is in fact

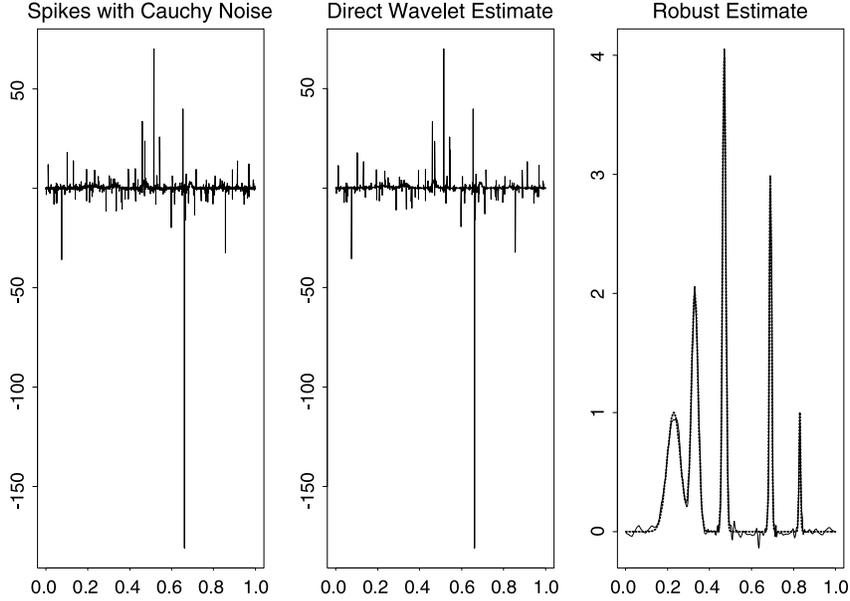

Fig. 1. *Left panel: spikes signal with Cauchy noise. Middle panel: an estimate obtained by applying directly a wavelet procedure to the original noisy signal. Right panel: a robust estimate by apply a wavelet block thresholding procedure to the medians of the binned data. Sample size is* 4096 *and bin size is* 8.

heavy-tailed. On the other hand, our robust procedure performs well even in Cauchy noise.

In Section 4 we construct a robust procedure for estimating the quadratic functional $Q(f) = \int f^2$ following the same general principles. Other problems such as construction of confidence sets and nonparametric hypothesis testing can be handled in a similar fashion.

Key technical tools used in our development are an improved moderate deviation result for the median statistic and a better quantile coupling inequality. Median coupling has been considered in Brown, Cai and Zhou (2008). For the asymptotic equivalence results given in Section 2 and the proofs of the theoretical results in Section 3 we need a more refined moderate deviation result for the median and an improved coupling inequality than those given in Brown, Cai and Zhou (2008). These improvements play a crucial role in this paper for establishing the asymptotic equivalence as well as robust and adaptive estimation results. The results may be of independent interest for other statistical applications.

The paper is organized as follows. Section 2 develops an asymptotic equivalence theory for unbounded loss functions. To illustrate the general principles of the asymptotic equivalence theory, we then consider robust estimation of the regression function $f$ under integrated squared error in Section 3 and estimation of the quadratic functional $\int f^2$ under squared error in Section 4. The two estimators are easily implementable and are shown to enjoy desirable robustness and adaptivity properties. In Section 5 we derive a moderate deviation result for the medians and a quantile coupling inequality. The proofs are contained in Section 6.

**2. Asymptotic equivalence.** This section develops an asymptotic equivalence theory for unbounded loss functions. The results reduce the general nonparametric regression model (1) to a standard Gaussian regression model.

The Gaussian nonparametric regression has been well studied and it often serves as a prototypical model for more general nonparametric function estimation settings. A large body of literature has been developed for minimax and adaptive estimation in the Gaussian case. These results include optimal convergence rates and optimal constants. See, for example, Pinsker (1980), Korostelev (1993), Donoho et al. (1995), Johnstone (2002), Tsybakov (2004), Cai and Low (2005, 2006b) and references therein for various estimation problems under various loss functions. The asymptotic equivalence results established in this section can be used to robustify these procedures in a unified way to treat the general nonparametric regression model (1).

We begin with a brief review of the classical formulation of asymptotic equivalence and then generalize it to accommodate unbounded losses.

2.1. *Classical asymptotic equivalence theory.* Le Cam (1986) developed a general theory for asymptotic decision problems. At the core of this theory is the concept of a distance between statistical models (or experiments), called Le Cam's deficiency distance. The goal is to approximate general statistical models by simple ones. If a complex model is close to a simple model in Le Cam's distance, then there is a mapping of solutions to decision theoretic problems from one model to the other for all *bounded* loss functions. Therefore the study of the complex model can be reduced to the one for the simple model.

A family of probability measures $E = \{P_\theta : \theta \in \Theta\}$ defined on the same $\sigma$-field of a sample space $\Omega$ is called a statistical model (or experiment). Le Cam (1964) defined a distance $\Delta(E, F)$ between $E$ and another model $F = \{Q_\theta : \theta \in \Theta\}$ with the same parameter set $\Theta$ by the means of "randomizations." Suppose one would like to approximate $E$ by a simpler model $F$. An observation $x$ in $E$ can be mapped into the sample space of $F$ by generating an "observation" $y$ according to a Markov kernel $K_x$, which is a probability measure on the sample space of $F$. Suppose $x$ is sampled from $P_\theta$. Write $KP_\theta$ for the distribution of $y$ with $KP_\theta(A) = \int K_x(A) \, dP_\theta$ for a measurable set $A$. The deficiency $\delta$ of $E$ with respect to $F$ is defined as the smallest possible value of the total variation distance between $KP_\theta$ and $Q_\theta$ among all possible choices of $K$, that is,

$$\delta(E,F) = \inf_K \sup_{\vartheta \in \Theta} |KP_\vartheta - Q_\vartheta|_{\mathrm{TV}}.$$

See Le Cam (1986, page 3) for further details. The deficiency $\delta$ of $E$ with respect to $F$ can be explained in terms of risk comparison. If $\delta(E,F) \leq \varepsilon$ for some $\varepsilon > 0$, it is easy to see that for every procedure $\tau$ in $F$ there exists a procedure $\xi$ in $E$ such that $R(\theta; \xi) \leq R(\theta; \tau) + 2\varepsilon$ for every $\theta \in \Theta$ and any loss function with values in the unit interval. The converse is also true. Symmetrically one may consider the deficiency of $F$ with respect to $E$ as well. The *Le Cam's deficiency distance* between the models $E$ and $F$ is then defined as

(3) $$\Delta(E,F) = \max(\delta(E,F), \delta(F,E)).$$

For bounded loss functions, if $\Delta(E,F)$ is small, then to every statistical procedure for $E$ there is a corresponding procedure for $F$ with almost the same risk function and vice versa. Two sequences of experiments $E_n$ and $F_n$ are called *asymptotically equivalent*, if $\Delta(E_n, F_n) \to 0$ as $n \to \infty$. The significance of asymptotic equivalence is that all asymptotically optimal statistical procedures can be carried over from one experiment to the other for bounded loss functions.

2.2. *Extension of the classical asymptotic equivalence formulation.* For many statistical applications, asymptotic equivalence under bounded losses is not sufficient because many commonly used loss functions are unbounded. Let $E_n = \{P_{\theta,n} : \theta \in \Theta\}$ and $F_n = \{Q_{\theta,n} : \theta \in \Theta\}$ be two asymptotically equivalent models in Le Cam's sense. Suppose that the model $F_n$ is simpler and well studied and a sequence of estimators $\hat\theta_n$ satisfy

$$\mathbb{E}_{Q_{\theta,n}} n^r d(\widehat\theta_n, \theta) \to c \qquad \text{as } n \to \infty,$$

where $d$ is a distance between $\widehat\theta$ and $\theta$, and $r, c > 0$ are constants. This implies that $\theta$ can be estimated by $\hat\theta_n$ under the distance $d$ with a rate $n^{-r}$. Examples include $\mathbb{E}_{Q_{\theta,n}} n(\widehat\theta - \theta)^2 \to c$ in many parametric estimation problems, and $\mathbb{E}_{Q_{f,n}} n^r \int (\widehat f - f)^2 \, d\mu \to c$, where $f$ is an unknown function and $0 < r < 1$, in many nonparametric estimation problems. The asymptotic equivalence between $E_n$ and $F_n$ in the classical sense does not imply that there is an estimator $\hat\theta^*$ in $E_n$ such that

$$\mathbb{E}_{P_{\theta,n}} n^r d(\hat\theta^*, \theta) \to c.$$

In this setting the loss function is actually $L(\vartheta, \theta) = n^r d(\vartheta, \theta)$ which grows as $n$ increases, and is usually unbounded.

In this section we introduce a new asymptotic equivalence formulation to handle unbounded losses. Let $\Lambda_E$ and $\Lambda_F$ be a set of procedures for $E$ and $F$, respectively. Let $\Gamma$ be a set of loss functions. We define the deficiency distance $\Delta(E, F; \Gamma, \Lambda_E, \Lambda_F)$ as follows.

DEFINITION 1. Define $\delta(E, F; \Gamma, \Lambda_E, \Lambda_F) \equiv \inf\{\epsilon \geq 0\colon$ for every procedure $\tau \in \Lambda_F$ there exists a procedure $\xi \in \Lambda_E$ such that $R(\theta; \xi) \leq R(\theta; \tau) + 2\varepsilon$ for every $\theta \in \Theta$ for any loss function $L \in \Gamma\}$. Then the deficiency distance between models $E$ and $F$ for the loss class $\Gamma$ and procedure classes $\Lambda_E$ and $\Lambda_F$ is defined as $\Delta(E, F; \Gamma, \Lambda_E, \Lambda_F) = \max\{\delta(E, F; \Gamma, \Lambda_E, \Lambda_F), \delta(F, E; \Gamma, \Lambda_F, \Lambda_E)\}$.

In other words, if the deficiency $\Delta(E, F; \Gamma, \Lambda_E, \Lambda_F)$ is small, then to every statistical procedure for one experiment, there is a corresponding procedure for another experiment with almost the same risk function for losses $L \in \Gamma$ and procedures in $\Lambda$.

DEFINITION 2. Two sequences of experiments $E_n$ and $F_n$ are called *asymptotically equivalent* with respect to the sets of procedures $\Lambda_{E_n}$ and $\Lambda_{F_n}$ and set of loss functions $\Gamma_n$ if $\Delta(E_n, F_n; \Gamma_n, \Lambda_{E_n}, \Lambda_{F_n}) \to 0$ as $n \to \infty$.

If $E_n$ and $F_n$ are asymptotically equivalent, then all asymptotically optimal statistical procedures in $\Lambda_{F_n}$ can be carried over to $E_n$ for loss functions $L \in \Gamma_n$ with essentially the same risk. The definitions here generalize the classical asymptotic equivalence formulation, which corresponds to the special case with $\Gamma$ being the set of loss functions with values in the unit interval.

For most statistical applications the loss function is bounded by a certain power of $n$. We now give a sufficient condition for the asymptotic equivalence under such losses. Suppose that we estimate $f$ or a functional of $f$ under a loss $L$. Let $p_{f,n}$ and $q_{f,n}$ be the density functions, respectively, for $E_n$ and $F_n$. Note that in the classical formulation of asymptotic equivalence for bounded losses, the deficiency of $E_n$ with respect to $F_n$ goes to zero if there is a Markov kernel $K$ such that

$$(4) \qquad \sup_f |KP_{f,n} - Q_{f,n}|_{\mathrm{TV}} \to 0.$$

For unbounded losses the condition (4) is no longer sufficient to guarantee that the deficiency goes to zero. Let $p^*_{f,n}$ and $q_{f,n}$ be the density functions of $KP_{f,n}$ and $Q_{f,n}$, respectively. Let $\varphi(f)$ be an estimand, which can be $f$ or a functional of $f$. Suppose that in $F_n$ there is an estimator $\widehat{\varphi(f)}_q$ of $\varphi(f)$ such that

$$\int L(\widehat{\varphi(f)}_q, \varphi(f)) q_{f,n} \to c.$$

We would like to derive sufficient conditions under which there is an estimator $\widehat{\varphi(f)}_p$ in $E_n$ such that

$$\int L(\widehat{\varphi(f)}_p, \varphi(f)) p_{f,n} \leq c(1 + o(1)).$$

Note that if $\widehat{\varphi(f)}_p$ is constructed by mapping over $\widehat{\varphi(f)}_q$ via a Markov kernel $K$, then

$$\mathbb{E} L(\widehat{\varphi(f)}_p, \varphi(f)) = \int L(\widehat{\varphi(f)}_q, \varphi(f)) p^*_{f,n}$$

$$= \int L(\widehat{\varphi(f)}_q, \varphi(f)) q_{f,n} + \int L(\widehat{\varphi(f)}_q, \varphi(f))(p^*_{f,n} - q_{f,n}).$$

Let $A_n = \{|p^*_{f,n}/q_{f,n} - 1| < \varepsilon_n\}$ for some $\varepsilon_n \to 0$, and write

$$\int L(\widehat{\varphi(f)}_q, \varphi(f))(p^*_{f,n} - q_{f,n})$$

$$= \int L(\widehat{\varphi(f)}_q, \varphi(f)) q_{f,n} (p^*_{f,n}/q_{f,n} - 1)[I(A_n) + I(A_n^c)]$$

$$\leq \int L(\widehat{\varphi(f)}_q, \varphi(f)) q_{f,n} (p^*_{f,n}/q_{f,n} - 1)\{p^*_{f,n}/q_{f,n} \geq 1\}[I(A_n) + I(A_n^c)]$$

$$\leq \varepsilon_n \int L(\widehat{\varphi(f)}_q, \varphi(f)) q_{f,n} I(A_n) + \int L(\widehat{\varphi(f)}_q, \varphi(f)) p^*_{f,n} I(A_n^c).$$

If $KP_{f,n}(A_n^c)$ decays exponentially fast uniformly over $\mathcal{F}$ and $L$ is bounded by a polynomial of $n$, this formula implies that

$$\int L(\widehat{\varphi(f)}_q, \varphi(f))(q_{f,n} - p^*_{f,n}) = o(1).$$

ASSUMPTION (A0). For each estimand $\varphi(f)$, each estimator $\widehat{\varphi(f)} \in \Lambda_n$ and each $L \in \Gamma_n$, there is a constant $M > 0$, independent of the loss function and the procedure, such that $L(\widehat{\varphi(f)}, \varphi(f)) \leq M n^M$.

The following result summarizes the above discussion and gives a sufficient condition for the asymptotic equivalence for the set of procedures $\Lambda_n$ and set of loss functions $\Gamma_n$.

PROPOSITION 1. *Let $E_n = \{P_{\theta,n} : \theta \in \Theta\}$ and $F_n = \{Q_{\theta,n} : \theta \in \Theta\}$ be two models. Suppose there is a Markov kernel $K$ such that $KP_{\theta,n}$ and $Q_{\theta,n}$ are defined on the same $\sigma$-field of a sample space. Let $p^*_{f,n}$ and $q_{f,n}$ be the density functions of $KP_{f,n}$ and $Q_{f,n}$ w.r.t. a dominating measure such that for a sequence $\varepsilon_n \to 0$*

$$\sup_f KP_{f,n}(|p^*_{f,n}/q_{f,n} - 1| \geq \varepsilon_n) \leq C_D n^{-D}$$

for all $D > 0$, then $\delta(E_n, F_n; \Gamma_n, \Lambda_{E_n}, \Lambda_{F_n}) \to 0$ as $n \to \infty$ under Assumption (A0).

Examples of loss functions include

$$L(\hat{f}_n, f) = n^{2\alpha/(2\alpha+1)} \int (\hat{f}_n - f)^2 \quad \text{and} \quad L(\hat{f}_n, f) = n^{2\alpha/(2\alpha+1)} \int (\sqrt{\hat{f}_n} - \sqrt{f})^2$$

for estimating $f$ and $L(\hat{f}_n, f) = n^{2\alpha/(2\alpha+1)} (\hat{f}_n(t_0) - f(t_0))^2$ for estimating $f$ at a fixed point $t_0$ where $\alpha$ is the smoothness of $f$, as long as we require $\hat{f}_n$ to be bounded by a power of $n$. If the maximum of $\hat{f}_n$ or $\hat{f}_n(t_0)$ grows faster than a polynomial of $n$, we commonly obtain a better estimate by truncation, for example, defining a new estimate $\min(\hat{f}_n, n^2)$.

The above discussions suggest that we may study a broad range of loss functions under a mild restriction on procedures. In comparison to the classic framework of asymptotic equivalence, here the collection of loss functions is much broadened to include unbounded losses while the collection of procedures is slightly more restrictive to only include those with losses bounded by a polynomial power of $n$. Virtually all practical procedures satisfy this condition. Of course in our formulation if the $\Gamma_n$ is set to be the collection of bounded loss functions, then the procedure can be any measurable function.

2.3. *Asymptotic equivalence for robust estimation under unbounded losses.* We now return to the nonparametric regression model (1) and denote the model by $E_n$,

$$E_n : Y_i = f(i/n) + \xi_i, \qquad i = 1, \ldots, n.$$

An asymptotic equivalence theory for nonparametric regression with a known error distribution has been developed in Grama and Nussbaum (2002), but the Markov kernel (randomization) there was not given explicitly, and so it is not implementable. In this section we propose an explicit and easily implementable procedure to reduce the nonparametric regression with an unknown error distribution to a Gaussian regression. We begin by dividing the interval $[0, 1]$ into $T$ equal-length subintervals. Without loss of generality, we shall assume that $n$ is divisible by $T$, and let $m = n/T$, the number of observations in each bin. We then take the median $X_j$ of the observations in each bin, that is,

$$X_j = \text{median}\{Y_i, (j-1)m + 1 \le i < jm\},$$

and make statistical inferences based on the median statistics $\{X_j\}$. Let $F_n$ be the experiment of observing $\{X_j, 1 \le j \le T\}$. In this section we shall show that $F_n$ is in fact asymptotically equivalent to the following Gaussian experiment:

$$G_n : X_j^{**} = f(j/T) + \frac{1}{2h(0)\sqrt{m}} Z_j, \qquad Z_j \overset{\text{i.i.d.}}{\sim} N(0,1),\ 1 \le j \le T,$$

under mild regularity conditions. The asymptotic equivalence is established in two steps.

Suppose the function $f$ is smooth. Then $f$ is locally approximately constant. We define a new experiment to approximate $E_n$ as follows:

$$E_n^*: Y_i^* = f^*(i/n) + \xi_i, \qquad 1 \leq i \leq n,$$

where $f^*(i/n) = f(\frac{[iT/n]}{T})$. For each of the $T$ subintervals, there are $m$ observations centered around the same mean.

For the experiment $E_n^*$ we bin the observations $Y_i^*$ and then take the medians in exactly the same way and let $X_j^*$ be the median of the $Y_i^*$'s in the $j$th subinterval. If $E_n^*$ approximates $E_n$ well, the statistical properties $X_j^*$ are then similar to $X_j$. Let $\eta_j$ be the median of corresponding errors $\xi_i$ in the $j$th bin. Note that the median of $X_j^*$ has a very simple form:

$$F_n^*: X_j^* = f(j/T) + \eta_j, \qquad 1 \leq j \leq T.$$

Theorem 6 in Section 5 shows that $\eta_j$ can be well approximated by a normal variable with mean 0 and variance $\frac{1}{4mh^2(0)}$, which suggests that $F_n^*$ is close to the experiment $G_n$.

We formalize the above heuristics in the following theorems. We first introduce some conditions. We shall choose $T = n^{2/3}/\log n$ and assume that $f$ is in a Hölder ball,

(5) $\qquad f \in \mathcal{F} = \{f : |f(y) - f(x)| \leq M|x-y|^d\}, \qquad d > 3/4.$

ASSUMPTION (A1). Let $\xi$ be a random variable with density function $h$. Define $r_a(\xi) = \log \frac{h(\xi-a)}{h(\xi)}$ and $\mu(a) = \mathbb{E}r(\xi)$. Assume that

(6) $\qquad\qquad\qquad\qquad \mu(a) \leq Ca^2,$

(7) $\qquad\qquad \mathbb{E}\exp[t(r_a(\xi) - \mu(a))] \leq \exp(Ct^2 a^2),$

for $0 \leq |a| < \varepsilon$ and $0 \leq |ta| < \varepsilon$ for some $\varepsilon > 0$. Equation (7) is roughly equivalent to $\mathbb{V}ar(r_a(\xi)) \leq Ca^2$. Assumption (A1) is satisfied by many distributions including Cauchy and Gaussian.

The following asymptotic equivalence result implies that any procedure based on $X_j$ has exactly the same asymptotic risk as a similar procedure by just replacing $X_j$ by $X_j^*$. That is, the experiments $F_n$ and $F_n^*$ are asymptotically equivalent.

THEOREM 1. *Under Assumptions (A0) and (A1) and the Hölder condition (5), the two experiments $E_n$ and $E_n^*$ are asymptotically equivalent with respect to the set of procedures $\Lambda_n$ and set of loss functions $\Gamma_n$.*

The following asymptotic equivalence result implies that asymptotically there is no need to distinguish $X_j^*$'s from the Gaussian random variables $X_j^{**}$'s. We need the following assumptions on the density function $h(x)$ of $\xi$.

ASSUMPTION (A2). $\int_{-\infty}^0 h(x) = \frac{1}{2}$, $h(0) > 0$, and $|h(x) - h(0)| \leq Cx^2$ in an open neighborhood of 0.

The last condition $|h(x) - h(0)| \leq Cx^2$ is basically equivalent to $h'(0) = 0$. The Assumption (A2) is satisfied when $h$ is symmetric and $h''$ exists in a neighborhood of 0.

THEOREM 2. *Under Assumptions (A0) and (A2), the two experiments $F_n^*$ and $G_n$ are asymptotically equivalent with respect to the set of procedures $\Lambda_n$ and set of loss functions $\Gamma_n$.*

These theorems together imply that, under assumptions (A1) and (A2) and the Hölder condition (5), the experiment $F_n$ is asymptotically equivalent to $G_n$ with respect to the set of procedures $\Lambda_n$ and set of loss functions $\Gamma_n$. So any statistical procedure $\delta$ in $G_n$ can be carried over to the $E_n$ (by treating $X_j$ as if it were $X_j^{**}$) in the sense that the new procedure has the same asymptotic risk as $\delta$ for all loss functions bounded by a certain power of $n$.

2.4. *Discussion.* The asymptotic equivalence theory provides deep insight and useful guidance for the construction of practical procedures in a broad range of statistical inference problems under the nonparametric regression model (1) with an unknown symmetric error distribution. Interesting problems include robust and adaptive estimation of the regression function, estimation of linear or quadratic functionals, construction of confidence sets, nonparametric hypothesis testing, etc. There is a large body of literature on these nonparametric problems in the case of Gaussian errors. With the asymptotic equivalence theory developed in this section, many of these procedures and results can be extended and robustified to deal with the case of an unknown symmetric error distribution. For example, the SureShrink procedure of Donoho and Johnstone (1995), the empirical Bayes procedures of Johnstone and Silverman (2005) and Zhang (2005), and SureBlock in Cai and Zhou (2009) can be carried over from the Gaussian regression to the general nonparametric regression. Theoretical properties such as rates of convergence remain the same under the regression model (1) with suitable regularity conditions.

To illustrate the general ideas, we consider in the next two sections two important nonparametric problems under the model (1): adaptive estimation of the regression function $f$ and robust estimation of the quadratic functional $Q(f) = \int f^2$. These examples show that for a given statistical problem it is easy to turn the case of nonparametric regression with general symmetric errors into the one with Gaussian noise and construct highly robust and adaptive procedures. Other robust inference problems can be handled in a similar fashion.

**3. Robust wavelet regression.** We consider in this section robust and adaptive estimation of the regression function $f$ under the model (1). Many estimation procedures have been developed in the literature for case where the errors $\xi_i$ are assumed to be i.i.d. Gaussian. However, these procedures are not readily applicable when the noise distribution is unknown. In fact direct application of the procedures designed for the Gaussian case can fail badly if the noise is in fact heavy-tailed. See, for example, Figure 1 in the Introduction.

In this section we construct a robust procedure by following the general principles of the asymptotic equivalence theory developed in Section 2. The estimator is robust, adaptive, and easily implementable. In particular, its performance is not sensitive to the error distribution.

3.1. *Wavelet procedure for robust nonparametric regression.* We begin with basic notation and definitions and then give a detailed description of our robust wavelet regression procedure.

Let $\{\phi, \psi\}$ be a pair of father and mother wavelets. The functions $\phi$ and $\psi$ are assumed to be compactly supported and $\int \phi = 1$. Dilation and translation of $\phi$ and $\psi$ generate an orthonormal wavelet basis. For simplicity in exposition, we work with periodized wavelet bases on $[0, 1]$. Let

$$\phi^p_{j,k}(t) = \sum_{l=-\infty}^{\infty} \phi_{j,k}(t-l), \qquad \psi^p_{j,k}(t) = \sum_{l=-\infty}^{\infty} \psi_{j,k}(t-l) \qquad \text{for } t \in [0,1],$$

where $\phi_{j,k}(t) = 2^{j/2}\phi(2^j t - k)$ and $\psi_{j,k}(t) = 2^{j/2}\psi(2^j t - k)$. The collection $\{\phi^p_{j_0,k}, k = 1, \ldots, 2^{j_0}; \psi^p_{j,k}, j \geq j_0 \geq 0, k = 1, \ldots, 2^j\}$ is then an orthonormal basis of $L^2[0,1]$, provided the primary resolution level $j_0$ is large enough to ensure that the support of the scaling functions and wavelets at level $j_0$ is not the whole of $[0, 1]$. The superscript "$p$" will be suppressed from the notation for convenience. An orthonormal wavelet basis has an associated orthogonal Discrete Wavelet Transform (DWT) which transforms sampled data into the wavelet coefficients. See Daubechies (1992) and Strang (1992) for further details on wavelets and discrete wavelet transform. A square-integrable function $f$ on $[0, 1]$ can be expanded into a wavelet series,

$$(8) \qquad f(t) = \sum_{k=1}^{2^{j_0}} \tilde{\theta}_{j_0,k} \phi_{j_0,k}(t) + \sum_{j=j_0}^{\infty} \sum_{k=1}^{2^j} \theta_{j,k} \psi_{j,k}(t),$$

where $\tilde{\theta}_{j,k} = \langle f, \phi_{j,k} \rangle, \theta_{j,k} = \langle f, \psi_{j,k} \rangle$ are the wavelet coefficients of $f$.

We now describe the robust regression procedure in detail. Let the sample $\{Y_i, i = 1, \ldots, n\}$ be given as in (1). Set $J = \lfloor \log_2 \frac{n}{\log^{1+b} n} \rfloor$ for some $b > 0$ and let $T = 2^J$. We first group the observations $Y_i$ consecutively into $T$ equilength bins and then take the median of each bin. Denote the medians by $X = (X_1, \ldots, X_T)$. Apply the discrete wavelet transform to the binned medians $X$ and let $U = T^{-1/2} W X$ be the empirical wavelet coefficients, where $W$ is the discrete wavelet transformation matrix. Write

$$(9) \quad U = (\tilde{y}_{j_0,1}, \ldots, \tilde{y}_{j_0,2^{j_0}}, y_{j_0,1}, \ldots, y_{j_0,2^{j_0}}, \ldots, y_{J-1,1}, \ldots, y_{J-1,2^{J-1}})'.$$

Here $\tilde{y}_{j_0,k}$ are the gross structure terms at the lowest resolution level, and $y_{j,k}$ ($j = j_0, \ldots, J-1, k = 1, \ldots, 2^j$) are empirical wavelet coefficients at level $j$ which represent fine structure at scale $2^j$. Set

$$(10) \quad \sigma_n = \frac{1}{2h(0)\sqrt{n}}.$$

Then the empirical wavelet coefficients can be written as

$$(11) \quad y_{j,k} = \theta_{j,k} + \epsilon_{j,k} + \sigma_n z_{j,k} + \xi_{j,k},$$

where $\theta_{j,k}$ are the true wavelet coefficients of $f$, $\epsilon_{j,k}$ are "small" deterministic approximation errors, $z_{j,k}$ are i.i.d. $N(0,1)$, and $\xi_{j,k}$ are some "small" stochastic errors. The asymptotic equivalence theory given in Section 2 indicates that both $\epsilon_{j,k}$ and $\xi_{j,k}$ are "negligible" and the calculations in Section 6 will show this is indeed the case. If these negligible errors are ignored then we have

$$(12) \quad y_{j,k} \approx \theta_{j,k} + \sigma_n z_{j,k} \qquad \text{with } z_{j,k} \overset{\text{i.i.d.}}{\sim} N(0,1),$$

which is the idealized Gaussian sequence model.

The BlockJS procedure introduced in Cai (1999) for Gaussian nonparametric regression is then applied to $y_{j,k}$ as if they are exactly distributed as in (12). More specifically, at each resolution level $j$, the empirical wavelet coefficients $y_{j,k}$ are grouped into nonoverlapping blocks of length $L$. Let $B_j^i = \{(j,k) : (i-1)L + 1 \le k \le iL\}$ and let $S_{j,i}^2 \equiv \sum_{(j,k) \in B_j^i} y_{j,k}^2$. Let $\hat{\sigma}_n^2$ be an estimator of $\sigma_n^2$ [see (16) for an estimator]. Set $J_* = \lfloor \log_2 \frac{T}{\log^{1+b} n} \rfloor$. A modified James–Stein shrinkage rule is then applied to each block $B_j^i$ with $j \le J_*$, that is,

$$(13) \quad \hat{\theta}_{j,k} = \left(1 - \frac{\lambda_* L \hat{\sigma}_n^2}{S_{j,i}^2}\right)_+ y_{j,k} \qquad \text{for } (j,k) \in B_j^i,$$

where $\lambda_* = 4.50524$ is a constant satisfying $\lambda_* - \log \lambda_* = 3$. For the gross structure terms at the lowest resolution level $j_0$, we set $\hat{\tilde{\theta}}_{j_0,k} = \tilde{y}_{j_0,k}$. The estimate of $f$ at the sample points $\{\frac{i}{T} : i = 1, \ldots, T\}$ is obtained by applying the inverse discrete wavelet transform (IDWT) to the denoised wavelet coefficients. That is, $\{f(\frac{i}{T}) : i = 1, \ldots, T\}$ is estimated by $\hat{f} = \{\widehat{f(\frac{i}{T})} : i = 1, \ldots, T\}$ with $\hat{f} = T^{1/2} W^{-1} \cdot \hat{\theta}$. The whole function $f$ is estimated by

$$\hat{f}_n(t) = \sum_{k=1}^{2^{j_0}} \hat{\tilde{\theta}}_{j_0,k} \phi_{j_0,k}(t) + \sum_{j=j_0}^{J_*-1} \sum_{k=1}^{2^j} \hat{\theta}_{j,k} \psi_{j,k}(t). \tag{14}$$

REMARK 1. An estimator of $h^{-2}(0)$ can be given by

$$\widehat{h}^{-2}(0) = \frac{8m}{T} \sum (X_{2k-1} - X_{2k})^2 \tag{15}$$

and the variance $\sigma_n^2$ is then estimated by

$$\hat{\sigma}_n^2 = \frac{1}{4\hat{h}^2(0)n} = \frac{2}{T^2} \sum (X_{2k-1} - X_{2k})^2. \tag{16}$$

It is shown in Section 6 that the estimator $\hat{\sigma}_n^2$ is an accurate estimate of $\sigma_n^2$.

The robust estimator $\hat{f}_n$ constructed above is easy to implement. Figure 2 below illustrate the main steps of the procedure. As a comparison, we also plotted the estimate obtained by applying directly the BlockJS procedure to the original noisy signal. It can be seen clearly that this wavelet procedure does not perform well in the case of heavy-tailed noise. Other standard wavelet procedures have similar performance qualitatively. On the other hand, the BlockJS procedure performs very well on the medians of the binned data.

3.2. *Adaptivity and robustness of the procedure.* The robust regression procedure presented in Section 3.1 enjoys a high degree of adaptivity and robustness. We consider the theoretical properties of the procedure over the Besov spaces. For a given $r$-regular mother wavelet $\psi$ with $r > \alpha$ and a fixed primary resolution level $j_0$, the Besov sequence norm $\|\cdot\|_{b_{p,q}^\alpha}$ of the wavelet coefficients of a function $f$ is defined by

$$\|f\|_{b_{p,q}^\alpha} = \|\underline{\xi}_{j_0}\|_p + \left(\sum_{j=j_0}^\infty (2^{js} \|\underline{\theta}_j\|_p)^q\right)^{1/q}, \tag{17}$$

where $\underline{\xi}_{j_0}$ is the vector of the father wavelet coefficients at the primary resolution level $j_0$, $\underline{\theta}_j$ is the vector of the wavelet coefficients at level $j$, and $s = \alpha + \frac{1}{2} - \frac{1}{p} > 0$. Note that the Besov function norm of index $(\alpha, p, q)$ of a function $f$ is equivalent to the sequence norm (17) of the wavelet coefficients of the function. See Meyer (1992), Triebel (1992) and DeVore and Popov (1988) for further details on Besov spaces. We define

$$B_{p,q}^{\alpha}(M) = \{f; \|f\|_{b_{p,q}^{\alpha}} \leq M\}. \tag{18}$$

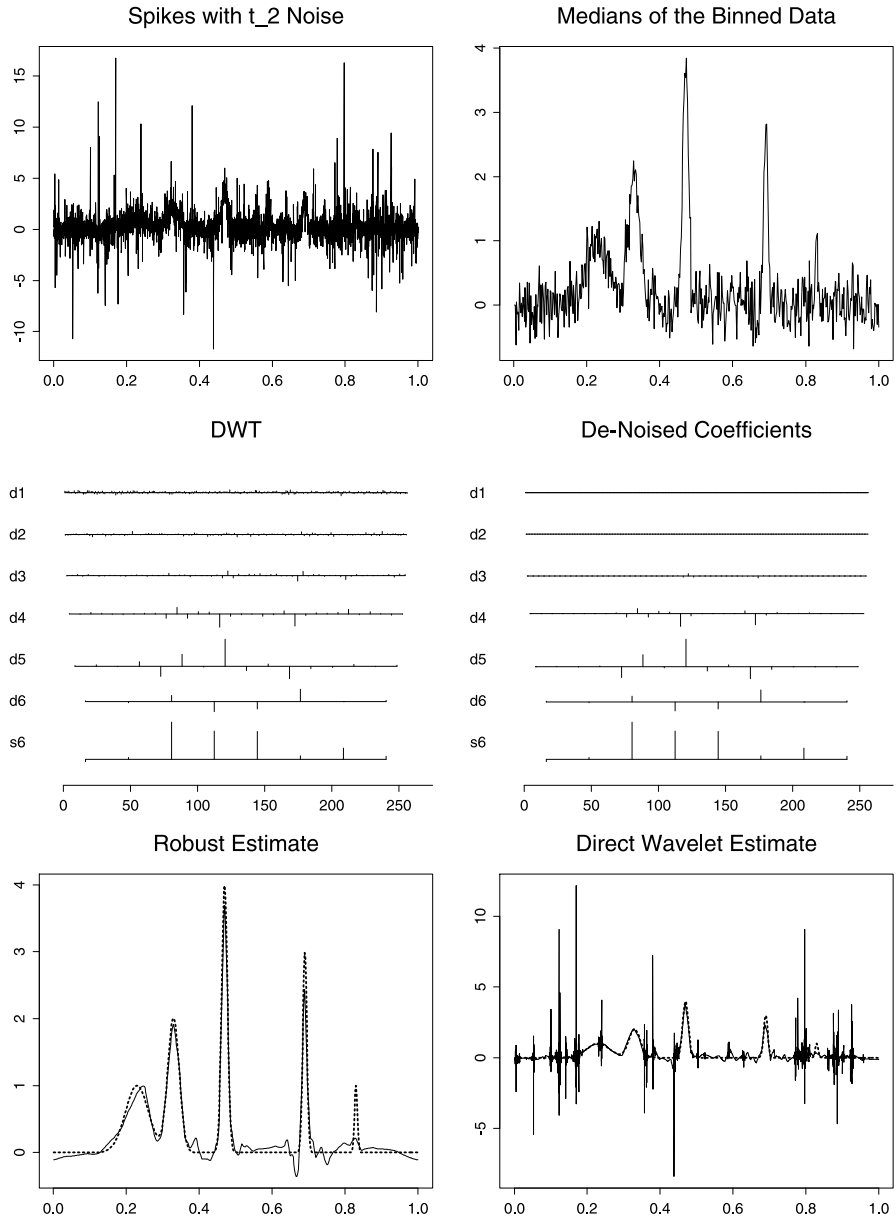

FIG. 2. *Top left panel: noisy Spikes signal with sample size $n = 4096$ where the noise has $t_2$ distribution. Top right panel: the medians of the binned data with the bin size $m = 8$. Middle left panel: the discrete wavelet coefficients of the medians. Middle right panel: blockwise thresholded wavelet coefficients of the medians. Bottom left panel: the robust estimate of the Spikes signal (dotted line is the true signal). Bottom right panel: the estimate obtained by applying directly the BlockJS procedure to the original noisy signal.*

In the case of Gaussian noise the minimax rate of convergence for estimating $f$ over the Besov body $B_{p,q}^{\alpha}(M)$ is $n^{-2\alpha/(1+2\alpha)}$. See Donoho and Johnstone (1998).

We shall consider the following collection of error distributions. For $0 < \epsilon_1 < 1$, $\epsilon_i > 0$, $i = 2, 3, 4$, let

$$\mathcal{H}_{\epsilon_1,\epsilon_2} = \left\{ h : \int_{-\infty}^{0} h(x) = \frac{1}{2}, \epsilon_1 \leq h(0) \leq \frac{1}{\epsilon_1}, \right.$$
$$\left. |h(x) - h(0)| \leq \frac{x^2}{\epsilon_1} \text{ for all } |x| < \epsilon_2 \right\} \tag{19}$$

and define $\mathcal{H} = \mathcal{H}(\epsilon_1, \epsilon_2, \epsilon_3, \epsilon_4)$ by

$$\mathcal{H} = \left\{ h \in \mathcal{H}_{\epsilon_1,\epsilon_2} : \int |x|^{\epsilon_3} h(x) \, dx < \epsilon_4, \right.$$
$$\left. h(x) = h(-x), |h^{(3)}(x)| \leq \epsilon_4 \text{ for } |x| \leq \epsilon_3 \right\}. \tag{20}$$

The assumption $\int |x|^{\epsilon_3} h(x) \, dx < \epsilon_4$ guarantees that the moments of the median of the binned data are well approximated by those of the normal random variable. Note that this assumption is satisfied by a large collection of distributions including Cauchy distribution.

The following theorem shows that our estimator achieves optimal global adaptation for a wide range of Besov balls $B_{p,q}^{\alpha}(M)$ defined in (18) and uniformly over the family of error distributions given in (20).

THEOREM 3. *Suppose the wavelet $\psi$ is $r$-regular. Then the estimator $\hat{f}_n$ defined in (14) satisfies, for $p \geq 2$, $\alpha \leq r$ and $\frac{2\alpha^2}{1+2\alpha} > \frac{1}{p}$,*

$$\sup_{h \in \mathcal{H}} \sup_{f \in B_{p,q}^{\alpha}(M)} E\|\hat{f}_n - f\|_2^2 \leq C n^{-2\alpha/(1+2\alpha)}$$

*and for $1 \leq p < 2$, $\alpha \leq r$ and $\frac{2\alpha^2}{1+2\alpha} > \frac{1}{p}$,*

$$\sup_{h \in \mathcal{H}} \sup_{f \in B_{p,q}^{\alpha}(M)} E\|\hat{f}_n - f\|_2^2 \leq C n^{-2\alpha/(1+2\alpha)} (\log n)^{(2-p)/(p(1+2\alpha))}.$$

In addition to global adaptivity, the estimator also enjoys a high degree of local spatial adaptivity. For a fixed point $t_0 \in (0, 1)$ and $0 < \alpha \leq 1$, define the local Hölder class $\Lambda^{\alpha}(M, t_0, \delta)$ by

$$\Lambda^{\alpha}(M, t_0, \delta) = \{f : |f(t) - f(t_0)| \leq M|t - t_0|^{\alpha} \text{ for } t \in (t_0 - \delta, t_0 + \delta)\}.$$

If $\alpha > 1$, then

$$\Lambda^\alpha(M, t_0, \delta) = \{f : |f^{(\lfloor\alpha\rfloor)}(t) - f^{(\lfloor\alpha\rfloor)}(t_0)| \leq M|t-t_0|^{\alpha'} \text{ for } t \in (t_0 - \delta, t_0 + \delta)\},$$

where $\lfloor\alpha\rfloor$ is the largest integer less than $\alpha$ and $\alpha' = \alpha - \lfloor\alpha\rfloor$.

In Gaussian nonparametric regression setting, it is well known that the optimal rate of convergence for estimating $f(t_0)$ over $\Lambda^\alpha(M, t_0, \delta)$ with $\alpha$ completely known is $n^{-2\alpha/(1+2\alpha)}$. On the other hand, when $\alpha$ is unknown, Lepski (1990) and Brown and Low (1996a) showed that the local adaptive minimax rate over the Hölder class $\Lambda^\alpha(M, t_0, \delta)$ is $(\log n/n)^{2\alpha/(1+2\alpha)}$. So one has to pay at least a logarithmic factor for adaptation.

Theorem 4 below shows that our estimator achieves optimal local adaptation with the minimal cost uniformly over the family of noise distributions defined in (20).

THEOREM 4. *Suppose the wavelet $\psi$ is $r$-regular with $r \geq \alpha > 0$. Let $t_0 \in (0,1)$ be fixed. Then the estimator $\hat{f}_n$ defined in (14) satisfies*

$$(21) \quad \sup_{h \in \mathcal{H}} \sup_{f \in \Lambda^\alpha(M, t_0, \delta)} E(\widehat{f}_n(t_0) - f(t_0))^2 \leq C \cdot \left(\frac{\log n}{n}\right)^{2\alpha/(1+2\alpha)}.$$

REMARK 2. Note that in the general asymptotic equivalence theory given in Section 2 the bin size was chosen to be $n^{1/3} \log n$. However, for specific estimation problems such as robust estimation of $f$ discussed in this section, the bin size can be chosen differently. Here we choose a small bin size $\log^{1+b} n$. There is a significant advantage in choosing such a small bin size in this problem. Note that the smoothness assumptions for $\alpha$ in Theorems 3 and 4 are different from those in Theorems 3 and 4 in Brown, Cai and Zhou (2008). For example, in Theorem 4 of Brown, Cai and Zhou (2008) it was assumed $\alpha > 1/6$, but now we need only $\alpha > 0$ due to the choice of the small bin size.

**4. Robust estimation of the quadratic functional $\int f^2$.** An important nonparametric estimation problem is that of estimating the quadratic functional $Q(f) = \int f^2$. This problem is interesting in its own right and closely related to the construction of confidence balls and nonparametric hypothesis testing in nonparametric function estimation. See, for example, Li (1989), Dümbgen (1998), Spokoiny (1998), Genovese and Wasserman (2005) and Cai and Low (2006a). In addition, as shown in Bickel and Ritov (1988), Donoho and Nussbaum (1990) and Fan (1991), this problem connects the nonparametric and semiparametric literatures.

Estimating the quadratic functional $Q(f)$ has been well studied in the Gaussian noise setting. See, for example, Donoho and Nussbaum (1990), Fan (1991). Efromovich and Low (1996), Laurent and Massart (2000), Klemelä (2006) and Cai and Low (2005, 2006b). In this section, we shall consider robust estimation of the quadratic functional $Q(f)$ under the regression model (1) with an unknown symmetric error distribution. We shall follow the same notation used in Section 3. Note that the orthonormality of the wavelet basis implies the isometry between the $L_2$ function norm and the $\ell_2$ wavelet sequence norm which yields

$$Q(f) = \int f^2 = \sum_{k=1}^{2^{j_0}} \tilde{\theta}_{j_0,k}^2 + \sum_{j=j_0}^{\infty} \sum_{k=1}^{2^j} \theta_{j,k}^2.$$

The problem of estimating $Q(f)$ is then translated into estimating the squared coefficients.

We consider adaptively estimating $Q(f)$ over Besov balls $B_{p,q}^\alpha(M)$ with $\alpha > \frac{1}{p} + \frac{1}{2}$. We shall show that it is in fact possible to find a simple procedure which is asymptotically rate optimal simultaneously over a large collection of unknown symmetric error distributions. In this sense, the procedure is robust.

As in Section 3, we group the observations $Y_i$ into $T$ bins of size $\log^{1+b}(n)$ for some $b > 0$ and then take the median of each bin. Let $X = (X_1, \ldots, X_T)$ denote the binned medians and let $U = T^{-1/2}WX$ be the empirical wavelet coefficients, where $W$ is the discrete wavelet transformation matrix. Write $U$ as in (9). Then the empirical wavelet coefficients can be approximately decomposed as in (12):

(22) $$\tilde{y}_{j_0,k} \approx \tilde{\theta}_{j_0,k} + \sigma_n \tilde{z}_{j_0,k} \quad \text{and} \quad y_{j,k} \approx \theta_{j,k} + \sigma_n z_{j,k},$$

where $\sigma_n = 1/(2h(0)\sqrt{n})$ and $\tilde{z}_{j_0,k}$ and $z_{j,k}$ are i.i.d. standard normal variables.

The quadratic functional $Q(f)$ can then be estimated as if we have exactly the idealized sequence model (22). More specifically, let $J_q = \lfloor \log_2 \sqrt{n} \rfloor$ and set

(23) $$\hat{Q} = \sum_{k=1}^{2^{j_0}} (\tilde{y}_{j_0,k}^2 - \hat{\sigma}_n^2) + \sum_{j=j_0}^{J_q} \sum_{k=1}^{2^j} (y_{j,k}^2 - \hat{\sigma}_n^2).$$

The following theorem shows that this estimator is robust and rate-optimal for a large collection of symmetric error distributions and a wide range of Besov classes simultaneously.

THEOREM 5. *For all Besov balls $B_{p,q}^\alpha(M)$ with $\alpha > \frac{1}{p} + \frac{1}{2}$, the estimator $\hat{Q}$ given in (23) satisfies*

(24) $$\sup_{f \in B_{p,q}^\alpha(M)} E_f(\hat{Q} - Q(f))^2 \leq \frac{M^2}{h^2(0)} n^{-1}(1 + o(1)).$$

REMARK 3. We should note that there is a slight tradeoff between efficiency and robustness. When the error distribution is known to be Gaussian, it is possible to construct a simple procedure which is efficient, asymptotically attaining the exact minimax risk $4M^2n^{-1}$. See, for example, Cai and Low (2005). In the Gaussian case, the upper bound in (24) is $2\pi M^2 n^{-1}$ which is slightly larger than $4M^2n^{-1}$. On the other hand, our procedure is robust over a large collection of unknown symmetric error distributions.

The examples of adaptive and robust estimation of the regression function and the quadratic functional given in the last and this sections illustrate the practical use of the general principles in the asymptotic equivalence theory given in Section 2. It is easy to see that other nonparametric inference problems such as the construction of confidence sets and nonparametric hypothesis testing under the general nonparametric regression model (1) can be handled in a similar way. Hence, our approach can be viewed as a general method for robust nonparametric inference.

**5. Technical tools: moderate deviation and quantile coupling for median.**
Quantile coupling is an important technical tool in probability and statistics. For example, the celebrated KMT coupling results given in Komlós, Major and Tusnády (1975) plays a key role in the Hungarian construction in the asymptotic equivalence theory. See, for example, Nussbaum (1996). Standard coupling inequalities are mostly focused on the coupling of the mean of i.i.d. random variables with a normal variable. Brown, Cai and Zhou (2008) studied the coupling of a median statistic with a normal variable. For the asymptotic equivalence theory given in Section 2 and the proofs of the theoretical results in Section 3 we need a more refined moderate deviation result for the median and an improved coupling inequality than those given in Brown, Cai and Zhou (2008). This improvement plays a crucial role in this paper. It is the main tool for reducing the problem of robust regression with unknown symmetric noise to a well studied and relatively simple problem of Gaussian regression. The result here may be of independent interest because of the fundamental role played by the median in statistics.

Let $X$ be a random variable with distribution $G$, and $Y$ with a continuous distribution $F$. Define

$$\widetilde{X} = G^{-1}(F(Y)), \tag{25}$$

where $G^{-1}(x) = \inf\{u : G(u) \geq x\}$, then $\mathcal{L}(\widetilde{X}) = \mathcal{L}(X)$. Note that $\widetilde{X}$ and $Y$ are now defined on the same probability space. This makes it possible to give a pointwise bound between $\widetilde{X}$ and $Y$. For example, one can couple Binomial$(m, 1/2)$ and $N(m/2, m/4)$ distributions. Let $X = 2(W - m/2)/\sqrt{m}$ with $W \sim \text{Binomial}(m, 1/2)$ and $Y \sim N(0, 1)$, and let $\widetilde{X}(Y)$ be defined as in (25). Komlós, Major and Tusnády (1975) showed that for some constant $C > 0$ and $\varepsilon > 0$, when $|\widetilde{X}| \leq \varepsilon\sqrt{m}$,

$$|\widetilde{X} - Y| \leq \frac{C}{\sqrt{m}} + \frac{C}{\sqrt{m}}|\widetilde{X}|^2. \tag{26}$$

Let $\xi_1, \ldots, \xi_m$ be i.i.d. random variables with density function $h$. Denote the sample median by $\xi_{\text{med}}$. The classical theory shows that the limiting distribution of $2h(0)\sqrt{m}\xi_{\text{med}}$ is $N(0,1)$. We will construct a new random variable $\widetilde{\xi}_{\text{med}}$ by using quantile coupling in (25) such that $\mathcal{L}(\widetilde{\xi}_{\text{med}}) = \mathcal{L}(\xi_{\text{med}})$ and show that $\widetilde{\xi}_{\text{med}}$ can be well approximated by a normal random variable as in (26). Denote the distribution and density function the sample median $\xi_{\text{med}}$ by $G$ and $g$, respectively. We obtain an improved approximation of the density $g$ by a normal density which leads to a better moderate deviation result for the distribution of sample median and consequently improve the classical KMT bound from the rate $1/\sqrt{m}$ to $1/m$. A general theory for improving the classical quantile coupling bound was given in Zhou (2006).

THEOREM 6. *Let $Z \sim N(0,1)$ and let $\xi_1, \ldots, \xi_m$ be i.i.d. with density function $h$, where $m = 2k + 1$ for some integer $k \geq 1$. Let Assumption (A2) hold. Then, for $|x| \leq \varepsilon$,*

$$g(x) = \frac{\sqrt{8k}f(0)}{\sqrt{2\pi}}\exp(-8kh^2(0)x^2/2 + O(kx^4 + k^{-1})) \tag{27}$$

*and for $0 < x < \varepsilon$,*

$$\begin{aligned}G(-x) &= \Phi(-x)\exp(O(kx^4 + k^{-1})) \quad and \\ \overline{G}(x) &= \overline{\Phi}(x)\exp(O(kx^4 + k^{-1})),\end{aligned} \tag{28}$$

*where $\overline{G}(x) = 1 - G(x)$, and $\overline{\Phi}(x) = 1 - \Phi(x)$. Consequently, for every $m$, there is a mapping $\widetilde{\xi}_{\text{med}}(Z): \mathbb{R} \mapsto \mathbb{R}$ such that $\mathcal{L}(\widetilde{\xi}_{\text{med}}(Z)) = \mathcal{L}(\xi_{\text{med}})$ and*

$$|2h(0)\sqrt{m}\widetilde{\xi}_{\text{med}} - Z| \leq \frac{C}{m} + \frac{C}{m}|2h(0)\sqrt{m}\widetilde{\xi}_{\text{med}}|^3, \qquad when \ |\widetilde{\xi}_{\text{med}}| \leq \varepsilon \tag{29}$$

*and*

$$|2h(0)\sqrt{m}\widetilde{\xi}_{\text{med}} - Z| \leq \frac{C}{m}(1 + |Z|^3), \qquad when \ |Z| \leq \varepsilon\sqrt{m}, \tag{30}$$

*where $C, \varepsilon > 0$ depend on $h$ but not on $m$.*

REMARK 4. In Brown, Cai and Zhou (2008), the density $g$ of the sample median was approximated by a normal density as

$$g(x) = \frac{\sqrt{8k}h(0)}{\sqrt{2\pi}}\exp(-8kh^2(0)x^2/2 + O(k|x|^3 + |x| + k^{-1})) \qquad for \ |x| \leq \varepsilon.$$

Since $\xi_{\text{med}} = O_p(1/\sqrt{m})$, the approximation error $O(k|x|^3 + |x| + k^{-1})$ is at the level of $1/\sqrt{m}$. In comparison, the approximation error $O(kx^4 + k^{-1})$ in (27) is at the level of $1/m$. This improvement is necessary for establishing (36) in the proof of Theorem 2, and leads to an improved quantile coupling bound (30) over the bound obtained in Brown, Cai and Zhou (2008):

$$|2h(0)\sqrt{m}\widetilde{\xi}_{\text{med}} - Z| \leq \frac{C}{\sqrt{m}} + \frac{C}{\sqrt{m}}Z^2, \qquad \text{when } |\widetilde{\xi}_{\text{med}}| \leq \varepsilon.$$

Since $Z$ is at a constant level, we improve the bound from a classical rate $1/\sqrt{m}$ to $1/m$.

Although the result is only given to $m$ odd, it can be easily extended to the even case as discussed in Remark 1 of Brown, Cai and Zhou (2008). The coupling result given in Theorem 6 in fact holds uniformly for the whole family of $h \in \mathcal{H}_{\epsilon_1,\epsilon_2}$.

THEOREM 7. *Let $\xi_1, \ldots, \xi_m$ be i.i.d. with density $h \in \mathcal{H}_{\epsilon_1,\epsilon_2}$ in (19). For every $m = 2k+1$ with integer $k \geq 1$, there is a mapping $\widetilde{\xi}_{\text{med}}(Z): \mathbb{R} \mapsto \mathbb{R}$ such that $\mathcal{L}(\widetilde{\xi}_{\text{med}}(Z)) = \mathcal{L}(\xi_{\text{med}})$ and for two constants $C_{\epsilon_1,\epsilon_2}, \varepsilon_{\epsilon_1,\epsilon_2} > 0$ depending only on $\epsilon_1$ and $\epsilon_2$,*

$$|2h(0)\sqrt{m}\widetilde{\xi}_{\text{med}} - Z| \leq \frac{C_{\epsilon_1,\epsilon_2}}{m} + \frac{C_{\epsilon_1,\epsilon_2}}{m}|2h(0)\sqrt{m}\widetilde{\xi}_{\text{med}}|^3, \tag{31}$$

$$\text{when } |\widetilde{\xi}_{\text{med}}| \leq \varepsilon_{\epsilon_1,\epsilon_2}$$

*and*

$$|2h(0)\sqrt{m}\widetilde{\xi}_{\text{med}} - Z| \leq \frac{C_{\epsilon_1,\epsilon_2}}{m} + \frac{C_{\epsilon_1,\epsilon_2}}{m}|Z|^3, \qquad \text{when } |Z| \leq \varepsilon_{\epsilon_1,\epsilon_2}\sqrt{m},$$

*uniformly over all $h \in \mathcal{H}_{\epsilon_1,\epsilon_2}$.*

**6. Proofs.** We shall prove the main results in the order of Theorems 6 and 7, Theorems 1 and 2, Theorem 3, and then Theorem 5. Theorems 6 and 7 provide important technical tools for the proof of the rest of the theorems. For reasons of space, we omit the proof of Theorem 4 and some of the technical lemmas. See Cai and Zhou (2008) for the complete proofs.

In this section, $C$ denotes a positive constant not depending on $n$ that may vary from place to place and we set $d \equiv \min(\alpha - \frac{1}{p}, 1)$.

6.1. *Proofs of Theorems 6 and 7.* We only prove (27) and (28). It follows from Zhou (2006) that the moderate deviation bound (28) implies the coupling bounds (29) and (30). Let $H(x)$ be the distribution function of $\xi_1$. The density of the median $\xi_{(k+1)}$ is

$$g(x) = \frac{(2k+1)!}{(k!)^2} H^k(x)(1-H(x))^k h(x).$$

Stirling's formula, $j! = \sqrt{2\pi} j^{j+1/2} \exp(-j + \epsilon_j)$ with $\epsilon_j = O(1/j)$, gives

$$g(x) = \frac{(2k+1)!}{4^k (k!)^2} [4H(x)(1-H(x))]^k h(x)$$

$$= \frac{2\sqrt{2k+1}}{e\sqrt{2\pi}} \left(\frac{2k+1}{2k}\right)^{2k+1} [4H(x)(1-H(x))]^k h(x) \exp\left(O\left(\frac{1}{k}\right)\right).$$

It is easy to see $|\sqrt{2k+1}/\sqrt{2k} - 1| \leq k^{-1}$, and

$$\left(\frac{2k+1}{2k}\right)^{2k+1} = \exp\left(-(2k+1)\log\left(1 - \frac{1}{2k+1}\right)\right) = \exp\left(1 + O\left(\frac{1}{k}\right)\right).$$

Then we have, when $0 < H(x) < 1$,

$$(32) \qquad g(x) = \frac{\sqrt{8k}}{\sqrt{2\pi}} [4H(x)(1-H(x))]^k h(x) \exp\left(O\left(\frac{1}{k}\right)\right).$$

From the assumption in the theorem, Taylor's expansion gives

$$4H(x)(1-H(x)) = 1 - 4(H(x) - H(0))^2$$

$$= 1 - 4\left[\int_0^x (h(t) - h(0))\,dt + h(0)x\right]^2$$

$$= 1 - 4(h(0)x + O(|x|^3))^2$$

for $0 \leq |x| \leq \varepsilon$, that is, $\log(4H(x)(1-H(x))) = -4h^2(0)x^2 + O(x^4)$ when $|x| \leq 2\varepsilon$ for some $\varepsilon > 0$. Here $\varepsilon$ is chosen sufficiently small so that $h(x) > 0$ for $|x| \leq 2\varepsilon$. Assumption (A2) also implies

$$\frac{h(x)}{h(0)} = 1 + O(x^2) = \exp(O(x^2)) \qquad \text{for } |x| \leq 2\varepsilon.$$

Thus, for $|x| \leq 2\varepsilon$,

$$g(x) = \frac{\sqrt{8k} h(0)}{\sqrt{2\pi}} \exp(-8kh^2(0)x^2/2 + O(kx^4 + x^2 + k^{-1}))$$

$$= \frac{\sqrt{8k} h(0)}{\sqrt{2\pi}} \exp(-8kh^2(0)x^2/2 + O(kx^4 + k^{-1})).$$

Now we approximate the distribution function of $\xi_{\text{med}}$ by a normal distribution. Without loss of generality, we assume $h(0) = 1$. We write

$$g(x) = \frac{\sqrt{8k}}{\sqrt{2\pi}} \exp(-8kx^2/2 + O(kx^4 + k^{-1})) \qquad \text{for } |x| \leq 2\varepsilon.$$

Now we use this approximation of density functions to give the desired approximation of distribution functions. Specifically, we shall show

$$(33) \qquad G(x) = \int_{-\infty}^{x} g(t)\,dt \le \Phi(\sqrt{8k}x)\exp(C(kx^4 + k^{-1}))$$

and

$$(34) \qquad G(x) \ge \Phi(\sqrt{8k}x)\exp(-C(kx^4 + k^{-1}))$$

for all $-\varepsilon \le x \le 0$ and some $C > 0$. The proof for $0 \le x \le \varepsilon$ is similar. We now prove inequality (33). Note that

$$(35) \qquad \begin{aligned} &(\Phi(\sqrt{8k}x)\exp(C(kx^4 + k^{-1})))' \\ &= \sqrt{8k}\varphi(\sqrt{8k}x)\exp(C(kx^4 + k^{-1})) \\ &\quad + \Phi(\sqrt{8k}x)4kCx^3\exp(C(kx^4 + k^{-1})). \end{aligned}$$

From Mill's ratio inequality, we have $\Phi(\sqrt{8k}x)(-\sqrt{8k}x) < \varphi(\sqrt{8k}x)$ and hence

$$\begin{aligned} &\Phi(\sqrt{8k}x)(4Ckx^3)\exp(C(kx^4 + k^{-1})) \\ &\ge \sqrt{8k}\varphi(\sqrt{8k}x)\left(-\frac{C}{2}x^2\right)\exp(C(kx^4 + k^{-1})). \end{aligned}$$

This and (35) yield

$$\begin{aligned} &(\Phi(\sqrt{8k}x)\exp(C(kx^4 + k^{-1})))' \\ &\ge \sqrt{8k}\varphi(\sqrt{8k}x)\left(1 - \frac{C}{2}x^2\right)\exp(C(kx^4 + k^{-1})) \\ &\ge \sqrt{8k}\varphi(\sqrt{8k}x)\exp(-Cx^2)\exp(C(kx^4 + k^{-1})) \\ &\ge \sqrt{8k}\varphi(\sqrt{8k}x)\exp(C(kx^4 + k^{-1})/4). \end{aligned}$$

Here in the second inequality we apply $1 - t/2 \ge \exp(-t)$ when $0 < t < 1/2$. Thus we have

$$(\Phi(\sqrt{8k}x)\exp(C(kx^4 + k^{-1})))' \ge \sqrt{8k}\varphi(\sqrt{8k}x)\exp(C(kx^4 + k^{-1}))$$

for $C$ sufficiently large and for $-2\varepsilon \le x \le 0$. Then

$$\begin{aligned} \int_{-2\varepsilon}^{x} g(t)\,dt &\le \int_{-2\varepsilon}^{x} (\Phi(\sqrt{8k}t)\exp(C(kx^4 + k^{-1})))' \\ &= \begin{bmatrix} \Phi(\sqrt{8k}x)\exp(C(kx^4 + k^{-1})) \\ -\Phi(\sqrt{8k}\cdot(2\varepsilon))\exp(C(k(2\varepsilon)^4 + k^{-1})) \end{bmatrix} \\ &\le \Phi(\sqrt{8k}x)\exp(C(kx^4 + k^{-1})). \end{aligned}$$

In (32) we see

$$\int_{-\infty}^{-2\varepsilon} g(t)\,dt = \int_{-\infty}^{-2\varepsilon} \frac{(2k+1)!}{(k!)^2} H^k(t)(1-H(t))^k h(t)\,dt$$

$$= \int_0^{H(-2\varepsilon)} \frac{(2k+1)!}{(k!)^2} u^k (1-u)^k\,du$$

$$= o(k^{-1}) \int_{H(-3\varepsilon/2)}^{H(-\varepsilon)} \frac{(2k+1)!}{(k!)^2} u^k (1-u)^k\,du$$

$$\leq o(k^{-1}) \int_{H(-2\varepsilon)}^{H(x)} \frac{(2k+1)!}{(k!)^2} u^k (1-u)^k\,du$$

$$= o(k^{-1}) \int_{-2\varepsilon}^{x} g(t)\,dt,$$

where the third equality is a result of the fact that $u_1^k(1-u_1)^k = o(k^{-1})u_2^k(1-u_2)^k$ uniformly for $u_1 \in [0, H(-2\varepsilon)]$ and $u_2 \in [H(-3\varepsilon/2), H(-\varepsilon)]$. Thus we have

$$G(x) \leq \Phi(\sqrt{8k}x) \exp(Ckx^4 + Ck^{-1}),$$

which is (33). Equation (34) can be established in a similar way.

REMARK. Note that in the proof of Theorem 6 it can be seen easily that constants $C$ and $\epsilon$ in (29) depends only on the ranges of $h(0)$ and the bound of Lipschitz constants of $h$ at a fixed open neighborhood of 0. Theorem 7 then follows from the proof of Theorem 6 together with this observation.

6.2. *Proofs of Theorems 1 and 2.*

PROOF OF THEOREM 1. Let $\varepsilon_n$ be a sequence approaching to 0 slowly, for example, $\varepsilon_n = 1/\log\log n$. Let $p_{f,n}$ be the joint density of $Y_i$'s and $p_{f,n}^*$ be the joint density of $Y_i^*$'s. And let $P_{f,n}$ be the joint distribution of $Y_i$'s and $P_{f^*,n}$ be the joint distribution of $Y_i^*$'s. We want to show that

$$\max\{P_{f^*,n}(|1 - p_{f^*,n}/p_{f,n}| \geq \varepsilon_n), P_{f,n}(|1 - p_{f,n}/p_{f^*,n}| \geq \varepsilon_n)\}$$

decays exponentially fast uniformly over the function space.

Note that $P_{f^*,n}(|1 - p_{f^*,n}/p_{f,n}| \geq \varepsilon_n) = P_{0,n}(|1 - p_{0,n}/p_{f^*-f,n}| \geq \varepsilon_n)$. It suffices to show that $P_{0,n}(|\log(p_{f^*-f,n}/p_{0,n})| \geq \varepsilon_n)$ decays exponentially fast. Write

$$\log(p_{f^*-f,n}/p_{0,n}) = \sum_{i=1}^n \log \frac{h(\xi_i - a_i)}{h(\xi_i)}$$

with $a_i = f^*(i/n) - f(i/n)$, where $\xi_i$ has density $h(x)$. Under Assumption (A1), we have $\mathbb{E} r_{a_i}(\xi_i) \leq C a_i^2$ and $\mathbb{E} \exp[t(r_{a_i}(\xi_i) - \mu(a_i))] \leq \exp(Ct^2 a_i^2)$ which imply

$$P_{0,n}\left(\exp\left[t\sum_{i=1}^n r_{a_i}(\xi_i) - \mu(a_i)\right] \geq \exp(t\varepsilon_n)\right) \leq \exp\left(Ct^2 \sum_{i=1}^n a_i^2 - t\varepsilon_n\right).$$

Since

$$\sum_{i=1}^n a_i^2 \leq C_1 n \cdot \left(\frac{n^{4/3}}{\log^2 n}\right)^{-d} = C_1 n^{1-4d/3} \log^{2d} n,$$

which goes to zero for $d > 3/4$, by setting $t = n^{(4d/3-1)/2}$ the Markov inequality implies that $P_{0,n}(|\log(p_{f^*-f,n}/p_{0,n})| \geq \varepsilon_n)$ decays exponentially fast. □

PROOF OF THEOREM 2. Let $g_{f,n}$ be the joint density of $X_j^*$'s and $q_{f,n}$ be the joint density of $X_j^{**}$'s. And let $G_{0,n}$ be the joint distribution of $\eta_j$'s and $Q_{0,n}$ be the joint distribution of $Z_j$'s. Theorem 6 yields

$$g(x) = \frac{\sqrt{4m}h(0)}{\sqrt{2\pi}} \exp(-4mh^2(0)x^2/2 + O(mx^4 + m^{-1}))$$

for $|x| \leq m^{-1/3}$. Since $G_{0,n}(|\eta_j| > m^{-1/3})$ and $Q_{0,n}(|Z_j| > m^{-1/3})$ decay exponentially fast, it suffices to study

$$\sum_{i=1}^T \log \frac{g(Z_j)}{\phi_{\sigma_m}(Z_j)} I(|Z_j| \leq m^{-1/3}).$$

Let

$$l(Z_j) = \log \frac{g(Z_j)}{\phi_{\sigma_m}(Z_j)} I(|Z_j| \leq m^{-1/3})$$

with $Z_j$ normally distributed with density $\phi_{\sigma_m}(x)$. It can be easily shown that

$$\mathbb{E} l(Z_j) \leq C Q_{0,n}\left(1 - \frac{g(Z_j)}{\phi_{\sigma_m}(Z_j)}\right)^2 \leq C_1 m^{-2}$$

and

$$\mathrm{Var}(l(Z_j)) \leq C m^{-2}.$$

Since $|Z_j| \leq Cm^{-1/3}$, then $|l(Z_j)| \leq Cm^{-1/3}$. Taylor's expansion gives

$$\mathbb{E} \exp[t(l(Z_j) - \mathbb{E} l(Z_j))] \leq \exp(Ct^2 m^{-2})$$

for $t = \log^{3/2} n$, then similar to the proof of Theorem 1 we have

(36) $\qquad Q_{f,n}(|\log(g_{f,n}/q_{f,n})| \geq \varepsilon_n) \leq \exp(Ct^2 T m^{-2} - t\varepsilon_n).$

Since $Tm^{-2} = 1/\log^3 n \to 0$, it decays faster than any polynomial of $n$. □

6.3. *Proof of Theorems 3 and 4.* In the proofs of Theorems 3, 4 and 5, we shall replace $\hat{\sigma}_n^2$ by $\sigma_n^2$. We assume that $h(0)$ is known and equal to 1 without loss of generality, since it can be shown easily that the estimator $\hat{h}(0)$ given in (15) satisfies

$$(37) \qquad P\{|\hat{h}^{-2}(0) - h^{-2}(0)| > n^{-\delta}\} \leq c_l n^{-l}$$

for some $\delta > 0$ and all constants $l \geq 1$. Note that $E\frac{m}{T/2}\sum(X_{2k-1} - X_{2k})^2 = \frac{1}{4}h^{-2}(0) + O(\sqrt{m}T^{-d})$, and it is easy to show

$$E\left|\frac{8m}{T}\sum(X_{2k-1} - X_{2k})^2 - h^{-2}(0)\right|^l \leq C_l(\sqrt{m}T^{-d})^l,$$

where $\sqrt{m}T^{-d} = n^{-\delta}$ with $\delta > 0$ in our assumption. Then (37) holds by Chebyshev's inequality. It is very important to see that the asymptotic risk properties of our estimators for $f$ in (13) and $Q(f)$ in (23) do not change when replacing $\sigma_n^2$ by $\sigma_n^2(1 + O(n^{-\delta}))$, thus in our analysis we may just assume that $h(0)$ is known without loss of generality.

For simplicity, we shall assume that $n$ is divisible by $T$ in the proof. The coupling inequality and the fact that a Besov ball $B_{p,q}^\alpha(M)$ can be embedded into a Hölder ball with smoothness $d = \min(\alpha - \frac{1}{p}, 1) > 0$ [see Meyer (1992)] enable us to precisely control of the errors. Proposition 2 gives the bounds for both the deterministic and stochastic errors.

PROPOSITION 2. *Let $X_j$ be given as in our procedure and let $f \in B_{p,q}^\alpha(M)$. Then $X_j$ can be written as*

$$(38) \qquad \sqrt{m}X_j = \sqrt{m}f\left(\frac{j}{T}\right) + \frac{1}{2}Z_j + \epsilon_j + \zeta_j,$$

*where:*

(i) $Z_j \overset{\text{i.i.d.}}{\sim} N(0, \frac{1}{h^2(0)})$;

(ii) $\epsilon_j$ *are constants satisfying* $|\epsilon_j| \leq C\sqrt{m}T^{-d}$ *and so* $\frac{1}{n}\sum_{i=1}^T \epsilon_j^2 \leq CT^{-2d}$;

(iii) $\zeta_j$ *are independent and "stochastically small" random variables satisfying with $E\zeta_j = 0$, and can be written as*

$$\zeta_j = \zeta_{j1} + \zeta_{j2} + \zeta_{j3}$$

*with*

$$|\zeta_{j1}| \leq C\sqrt{m}T^{-d},$$

$$E\zeta_{j2} = 0 \quad and \quad |\zeta_{2j}| \leq \frac{C}{m}(1 + |Z_j|^3),$$

$$P(\zeta_{j3} = 0) \geq 1 - C\exp(-\varepsilon m) \quad and \quad E|\zeta_{j3}|^D \text{ exists}$$

*for some $\varepsilon > 0$ and $C > 0$, and all $D > 0$.*

REMARK 5. Equation (38) is different than Proposition 1 in Brown, Cai and Zhou (2008), where there is an additional bias term $\sqrt{m}b_m$. Lemma 5 in Brown, Cai and Zhou (2008) showed that the bias $b_m$ can be estimated with a rate $\max\{T^{-2d}, m^{-4}\}$. Therefore in that paper we need to choose the bin size $m = n^{1/4}$ such that $m^{-4} = o(n^{-2\alpha/(2\alpha+1)})$ is negligible relative to the minimax risk. In the present paper we can choose $m = \log^{1+b} n$ because there is no bias term and as a result the condition on the smoothness is relaxed.

The proof of Proposition 2 is similar to that of Proposition 1 in Brown, Cai and Zhou (2008) and is thus omitted here. See Cai and Zhou (2008) for a complete proof.

We now consider the wavelet transform of the medians of the binned data. From Proposition 2 we may write

$$\frac{1}{\sqrt{T}}X_i = \frac{f(i/T)}{\sqrt{T}} + \frac{\epsilon_i}{\sqrt{n}} + \frac{Z_i}{2\sqrt{n}} + \frac{\zeta_i}{\sqrt{n}}.$$

Let $(y_{j,k}) = T^{-1/2}W \cdot X$ be the discrete wavelet transform of the binned data. Then one may write

$$(39) \qquad y_{j,k} = \theta'_{j,k} + \epsilon_{j,k} + \frac{1}{2\sqrt{n}}z_{j,k} + \xi_{j,k},$$

where $\theta'_{j,k}$ are the discrete wavelet transform of $(f(\frac{i}{T}))_{1\leq i\leq T}$, $z_{j,k}$ are the transform of the $Z_i$'s and so are i.i.d. $N(0,1)$ and $\epsilon_{j,k}$ and $\xi_{j,k}$ are, respectively, the transforms of $(\frac{\epsilon_i}{\sqrt{n}})$ and $(\frac{\zeta_i}{\sqrt{n}})$. The following proposition gives the risk bounds of the block thresholding estimator in a single block. These risk bounds are similar to results for the Gaussian case given in Cai (1999). But in the current setting the error terms $\epsilon_{j,k}$ and $\xi_{j,k}$ make the problem more complicated.

PROPOSITION 3. *Let $y_{j,k}$ be given as in (39) and let the block thresholding estimator $\hat{\theta}_{j,k}$ be defined as in (13). Then:*

(i) *for some constant $C > 0$,*

$$(40) \qquad E\sum_{(j,k)\in B_j^i}(\hat{\theta}_{j,k} - \theta'_{j,k})^2 \leq \min\left\{4\sum_{(j,k)\in B_j^i}(\theta'_{j,k})^2, 8\lambda_* L n^{-1}\right\} + 6\sum_{(j,k)\in B_j^i}\epsilon_{j,k}^2 + CLn^{-2};$$

(ii) *for any $0 < \tau < 1$, there exists a constant $C_\tau > 0$ depending on $\tau$ only such that for all $(j,k) \in B_j^i$*

(41) $\quad E(\hat{\theta}_{j,k} - \theta'_{j,k})^2 \leq C_\tau \cdot \min\Big\{\max_{(j,k)\in B_j^i}\{(\theta'_{j,k} + \epsilon_{j,k})^2\}, Ln^{-1}\Big\} + n^{-2+\tau};$

(iii) *for $j \leq J_*$ and $\epsilon_n > 1/\log n$, $P(\sqrt{n}|\xi_{j,k}| \geq \varepsilon_n) \leq C\exp(-\varepsilon_n m)$.*

The third part follows from Lemma 3 in Cai and Wang (2008) which gives a concentration inequality for wavelet coefficients at a given resolution.

For reasons of space we omit the proof of Proposition 3 here. See Cai and Zhou (2008) for a complete proof. We also need the following lemmas for the proof of Theorems 3 and 4. The proof of these lemmas is relatively straightforward and is thus omitted.

LEMMA 1. *Suppose $y_i = \theta_i + z_i, i = 1, \ldots, L$, where $\theta_i$ are constants and $z_i$ are random variables. Let $S^2 = \sum_{i=1}^L y_i^2$ and let $\hat{\theta}_i = (1 - \frac{\lambda L}{S^2})_+ y_i$. Then*

(42) $\quad E\|\hat{\theta} - \theta\|_2^2 \leq \|\theta\|_2^2 \wedge 4\lambda L + 4E[\|z\|_2^2 I(\|z\|_2^2 > \lambda L)].$

LEMMA 2. *Let $X \sim \chi_L^2$ and $\lambda > 1$. Then*

(43) $$P(X \geq \lambda L) \leq e^{-L/2(\lambda - \log \lambda - 1)} \quad \text{and}$$
$$EXI(X \geq \lambda L) \leq \lambda L e^{-L/2(\lambda - \log \lambda - 1)}.$$

LEMMA 3. *Let $T = 2^J$ and let $f_J(x) = \sum_{k=1}^T \frac{1}{\sqrt{T}} f(\frac{k}{T}) \phi_{J,k}(x)$. Then*

$$\sup_{f \in B_{p,q}^\alpha(M)} \|f_J - f\|_2^2 \leq CT^{-2d} \quad \text{where } d = \min(\alpha - 1/p, 1).$$

*Let $\{\theta'_{j,k}\}$ be the discrete wavelet transform of $\{f(\frac{i}{T}), 1 \leq i \leq T\}$ and let $\{\theta_{j,k}\}$ be the true wavelet coefficients of $f$. Then $|\theta'_{j,k} - \theta_{j,k}| \leq CT^{-d} 2^{-j/2}$ and consequently $\sum_{j=j_0}^{J-1} \sum_k (\theta'_{j,k} - \theta_{j,k})^2 \leq CT^{-2d}$.*

6.3.1. *Global adaptation: proof of Theorem 3.* Decompose $E\|\hat{f}_n - f\|_2^2$ into three terms as follows:

$$E\|\hat{g}_n - g\|_2^2 = \sum_k E(\hat{\tilde{\theta}}_{j_0,k} - \tilde{\theta}_{j,k})^2 + \sum_{j=j_0}^{J_*-1} \sum_k E(\hat{\theta}_{j,k} - \theta_{j,k})^2$$

(44) $$+ \sum_{j=J_*}^\infty \sum_k \theta_{j,k}^2$$

$$\equiv S_1 + S_2 + S_3.$$

It is easy to see that the first term $S_1$ and the third term $S_3$ are small:

$$S_1 = 2^{j_0} n^{-1} \epsilon^2 = o(n^{-2\alpha/(1+2\alpha)}). \tag{45}$$

Note that for $x \in \mathbb{R}^m$ and $0 < p_1 \leq p_2 \leq \infty$,

$$\|x\|_{p_2} \leq \|x\|_{p_1} \leq m^{1/p_1 - 1/p_2} \|x\|_{p_2}. \tag{46}$$

Since $f \in B^\alpha_{p,q}(M)$, so $2^{js} (\sum_{k=1}^{2^j} |\theta_{j,k}|^p)^{1/p} \leq M$. Now (46) yields that

$$S_3 = \sum_{j=J_*}^{\infty} \sum_k \theta_{j,k}^2 \leq C 2^{-2J_*(\alpha \wedge (\alpha + 1/2 - 1/p))}. \tag{47}$$

Propositions 2(ii) and 3 and Lemma 3 together yield

$$\begin{aligned}
S_2 &\leq 2 \sum_{j=j_0}^{J_*-1} \sum_k E(\hat{\theta}_{j,k} - \theta'_{j,k})^2 + 2 \sum_{j=j_0}^{J_*-1} \sum_k (\theta'_{j,k} - \theta_{j,k})^2 \\
&\leq \sum_{j=j_0}^{J_*-1} \sum_{i=1}^{2^j/L} \min\Big\{8 \sum_{(j,k) \in B_j^i} \theta_{j,k}^2, 8\lambda_* L n^{-1}\Big\} + 6 \sum_{j=j_0}^{J_*-1} \sum_k \epsilon_{j,k}^2 \\
&\quad + C n^{-1} + 10 \sum_{j=j_0}^{J_*-1} \sum_k (\theta'_{j,k} - \theta_{j,k})^2 \\
&\leq \sum_{j=j_0}^{J_*-1} \sum_{i=1}^{2^j/L} \min\Big\{8 \sum_{(j,k) \in B_j^i} \theta_{j,k}^2, 8\lambda_* L n^{-1}\Big\} + C n^{-1} + C T^{-2d}.
\end{aligned} \tag{48}$$

We now divide into two cases. First consider the case $p \geq 2$. Let $J_1 = [\frac{1}{1+2\alpha} \times \log_2 n]$. So, $2^{J_1} \approx n^{1/(1+2\alpha)}$. Then (48) and (46) yield

$$\begin{aligned}
S_2 &\leq 8\lambda_* \sum_{j=j_0}^{J_1-1} \sum_{i=1}^{2^j/L} L n^{-1} + 8 \sum_{j=J_1}^{J_*-1} \sum_k \theta_{j,k}^2 + C n^{-1} + C T^{-2d} \\
&\leq C n^{-2\alpha/(1+2\alpha)}.
\end{aligned}$$

By combining this with (45) and (47), we have $E\|\hat{f}_n - f\|_2^2 \leq C n^{-2\alpha/(1+2\alpha)}$ for $p \geq 2$.

Now let us consider the case $p < 2$. First we state the following lemma without proof.

LEMMA 4. *Let $0 < p < 1$ and $S = \{x \in \mathbb{R}^k : \sum_{i=1}^k x_i^p \leq B, x_i \geq 0, i = 1, \ldots, k\}$. Then for $A > 0$, $\sup_{x \in S} \sum_{i=1}^k (x_i \wedge A) \leq B \cdot A^{1-p}$.*

Let $J_2$ be an integer satisfying $2^{J_2} \asymp n^{1/(1+2\alpha)}(\log n)^{(2-p)/p(1+2\alpha)}$. Note that

$$\sum_{i=1}^{2^j/L}\left(\sum_{(j,k)\in B_j^i}\theta_{j,k}^2\right)^{p/2} \leq \sum_{k=1}^{2^j}(\theta_{j,k}^2)^{p/2} \leq M2^{-jsp}.$$

It then follows from Lemma 4 that

(49)
$$\sum_{j=J_2}^{J_*-1}\sum_{i=1}^{2^j/L}\min\left\{8\sum_{(j,k)\in B_j^i}\theta_{j,k}^2, 8\lambda_*Ln^{-1}\right\}$$
$$\leq Cn^{-2\alpha/(1+2\alpha)}(\log n)^{(2-p)/(p(1+2\alpha))}.$$

On the other hand,

(50)
$$\sum_{j=j_0}^{J_2-1}\sum_{i=1}^{2^j/L}\min\left\{8\sum_{(j,k)\in B_j^i}\theta_{j,k}^2, 8\lambda_*Ln^{-1}\right\}$$
$$\leq \sum_{j=j_0}^{J_2-1}\sum_b 8\lambda_*Ln^{-1}$$
$$\leq Cn^{-2\alpha/(1+2\alpha)}(\log n)^{(2-p)/(p(1+2\alpha))}.$$

We finish the proof for the case $p < 2$ by putting (45), (47), (49) and (50) together:

$$E\|\hat{f}_n - f\|_2^2 \leq Cn^{-2\alpha/(1+2\alpha)}(\log n)^{(2-p)/(p(1+2\alpha))}.$$

6.4. *Proof of Theorem 5.* Recall that

$$\hat{Q} = \sum_{k=1}^{2^{j_0}}(\tilde{y}_{j_0,k}^2 - \hat{\sigma}_n^2) + \sum_{j=j_0}^{J_q}\sum_{k=1}^{2^j}(y_{j,k}^2 - \hat{\sigma}_n^2)$$

and note that the empirical wavelet coefficients can be written as

$$y_{j,k} = \theta_{j,k} + \epsilon_{j,k} + \sigma_n z_{j,k} + \xi_{j,k}.$$

Since $(\sum_{j>J_1}\theta_{j,k}^2)^2 \leq C[2^{-2J_1(\alpha-1/p)}]^2 = o(\frac{1}{n})$, as in Cai and Low (2005) it is easy to show that

$$E_f\left\{\sum_{k=1}^{2^{j_0}}[(\tilde{\theta}_{j_0,k} + \sigma_n\tilde{z}_{j_0,k})^2 - \sigma_n^2] + \sum_{j=j_0}^{J_q}\sum_{k=1}^{2^j}[(\theta_{j,k} + \sigma_n z_{j,k})^2 - \sigma_n^2] - Q(f)\right\}^2$$
$$\leq 4\sigma_n^2 M^2(1 + o(1)).$$

The theorem then follows easily from the facts below:

$$\sum_{k=1}^{2^{j_0}} \tilde{\epsilon}_{j,k}^2 + \sum_{j=j_0}^{J_q} \sum_{k=1}^{2^j} \epsilon_{j,k}^2 \leq CT^{-2(\alpha-1/p)} = o\left(\frac{1}{n}\right),$$

$$E\left(\sum_{k=1}^{2^{j_0}} \xi_{j,k}^2 + \sum_{j=j_0}^{J_q} \sum_{k=1}^{2^j} \xi_{j,k}^2\right)^2 \leq C \frac{1}{m^2 n} = o\left(\frac{1}{n}\right),$$

$$E[\sqrt{n}(\hat{\sigma}_n^2 - \sigma_n^2)]^2 = o\left(\frac{1}{n}\right).$$

Department of Statistics  
The Wharton School  
University of Pennsylvania  
Philadelphia, Pennsylvania 19104  
USA  
E-mail: tcai@wharton.upenn.edu

Department of Statistics  
P.O. Box 208290  
Yale University  
New Haven, Connecticut 06511  
USA  
E-mail: huibin.zhou@yale.edu